\documentclass{llncs}

\usepackage{amsmath,amssymb,amsfonts}
\begin{document}

\title{The inverse of the star-discrepancy problem and the generation of pseudo-random numbers}

\author{Josef Dick and Friedrich Pillichshammer}

\institute{School of Mathematics and Statistics, The University of New South Wales, Sydney NSW Australia; email: josef.dick@unsw.edu.au \\ and \\ Department of Financial Mathematics, Johannes Kepler University, Linz, Austria; email: friedrich.pillichshammer@jku.at}

\maketitle

\begin{abstract}
The inverse of the star-discrepancy problem asks for point sets $P_{N,s}$ of size $N$ in the $s$-dimensional unit cube $[0,1]^s$ whose star-discrepancy $D^\ast(P_{N,s})$ satisfies $$D^\ast(P_{N,s}) \le C \sqrt{s/N},$$ where $C> 0$ is a constant independent of $N$ and $s$. The first existence results in this direction were shown by Heinrich, Novak, Wasilkowski, and Wo\'{z}niakowski in 2001, and a number of improvements have been shown since then. Until now only proofs that such point sets exist are known. Since such point sets would be useful in applications, the big open problem is to find explicit constructions of suitable point sets $P_{N,s}$.

We review the current state of the art on this problem and point out some connections to pseudo-random number generators.
\end{abstract}

\section{Introduction}\label{sec:Introduction}

The star-discrepancy is a quantitative measure for the irregularity of distribution of a point set $P_{N,s} = \{ \boldsymbol{x}_0, \boldsymbol{x}_1, \ldots, \boldsymbol{x}_{N-1}\}$ in the $s$-dimensional unit cube $[0,1)^s$. It is defined as the $L_\infty$-norm of the local discrepancy $$\Delta(\boldsymbol{t}) := \frac{\#\{n \in \{0,1,\ldots,N-1\} \ : \ \boldsymbol{x}_n \in [\boldsymbol{0},\boldsymbol{t})\}}{N}- \lambda_s([\boldsymbol{0},\boldsymbol{t})),$$ for $\boldsymbol{t}=(t_1,t_2,\ldots,t_s) \in [0,1]^s$, where $[\boldsymbol{0},\boldsymbol{t})=\prod_{j=1}^s[0,t_j)$ and $\lambda_s$ denotes the $s$-dimensional Lebesgue measure. In other words, the star-discrepancy (or $L_\infty$-discrepancy) of $P_{N,s}$ is $$D_N^{\ast}(P_{N,s})=\sup_{\boldsymbol{t} \in [0,1]^s} |\Delta(\boldsymbol{t})|.$$ Its significance arises from the classical Koksma-Hlawka inequality \cite{Hlawka,Koksma}, which states that
\begin{equation}\label{Koksma-Hlawka_inequality}
\left|\int_{[0,1]^s} f(\boldsymbol{x}) \,\mathrm{d} \boldsymbol{x} - \frac{1}{N} \sum_{n=0}^{N-1} f(\boldsymbol{x}_n) \right| \le V(f) D^\ast(P_{N,s}),
\end{equation}
where $V(f)$ denotes the variation of $f$ in the sense of Hardy and Krause, see, e.g., \cite{DP10,kuinie,niesiam}. This is the fundamental error estimate for quasi-Monte Carlo rules $Q(f)=(1/N) \sum_{n=0}^{N-1} f(\boldsymbol{x}_n)$.

To provide some insight into this inequality, we prove a simple version of it. Let $f:[0,1] \to \mathbb{R}$ be absolutely continuous, then for $x \in [0,1]$ we have
\begin{equation}\label{f_rep}
f(x) = f(1) - \int_0^1 \boldsymbol{1}_{[x,1]}(t) f'(t) \,\mathrm{d} t,
\end{equation}
where $\boldsymbol{1}$ denotes the indicator function. Using \eqref{f_rep} we obtain
\begin{align*}
\int_{[0,1]^s} f(x) \,\mathrm{d} x - \frac{1}{N} \sum_{n=0}^{N-1} f(x_n)  = & \int_0^1 f'(t) \left[ \frac{1}{N} \sum_{n=0}^{N-1} \boldsymbol{1}_{[x_n,1]}(t) -  \int_0^1 \boldsymbol{1}_{[x,1]}(t) \,\mathrm{d} x \right] \,\mathrm{d} t \\ = & \int_0^1 f'(t) \left[\frac{1}{N} \sum_{n=0}^{N-1} \boldsymbol{1}_{[0, t]}(x_n) - t \right] \,\mathrm{d} t.
\end{align*}
This implies that
\begin{equation*}
\left| \int_{[0,1]^s} f(x) \,\mathrm{d} x - \frac{1}{N} \sum_{n=0}^{N-1} f(x_n)  \right| \le \int_0^1 |f'(t)| \,\mathrm{d} t \sup_{0 \le t \le 1} \left|\frac{1}{N} \sum_{n=0}^{N-1} \boldsymbol{1}_{[0,t]}(x_n) - t \right|.
\end{equation*}
The right-most expression in the above inequality is simply the star-discrepancy of the point set $\{x_0, x_1,\ldots, x_{N-1}\}$ and for absolutely continuous functions $f$, the term $\int_0^1 |f'(t)| \,\mathrm{d} t$ coincides with the Hardy-Krause variation. This approach can be generalized to the $s$-dimensional unit cube $[0,1]^s$, yielding a version of the Koksma-Hlawka inequality \eqref{Koksma-Hlawka_inequality}. To obtain quasi-Monte Carlo rules with small quadrature error, it is therefore of importance to design point sets with small star-discrepancy.

In many papers the star-discrepancy is studied from the view point of its asymptotic behavior in $N$ (for a fixed dimension $s$). Define the {\it $N$th minimal star-discrepancy} in $[0,1)^s$ as $${\rm disc}(N,s):= \inf_{P_{N,s}} D_N^{\ast}(P_{N,s}),$$ where the infimum is extended over all $N$-element point sets in $[0,1)^s$. It is well known that ${\rm disc}(N,s)$ behaves like 
\begin{equation}\label{star_bounds}
\frac{(\log N)^{(s-1)/2 + \delta_s}}{N} \ll_s {\rm disc}(N,s) \ll_s \frac{(\log N)^{s-1}}{N}, 
\end{equation}
where $\delta_s \in (0,1/2)$ is an unknown quantity depending only on $s$. Here $A \ll_s B$ means that there is a constant $c_s>0$ depending only on $s$ such that $A \le c_s B$. The lower bound was shown by Bilyk, Lacey and Vagharshakyan \cite{blv08} improving a famous result of Roth~\cite{roth1}. For the upper bound several explicit constructions are known whose star-discrepancy achieves such a bound. See, e.g., \cite{DP10,niesiam}. Thus the upper bound on the $N$th minimal star-discrepancy is of order of magnitude $O(N^{-1+\varepsilon})$ for every $\varepsilon>0$. The problem however is that the function $N \mapsto (\log N)^{s-1}/N$ does not decreases to zero until $N> \exp(s-1)$. For $N \le \exp(s-1)$ this function is increasing which means that for $N$ in this range our discrepancy bound is useless. Even for moderately large dimensions $s$, point sets with cardinality $N > \exp(s-1)$ cannot be used for practical applications in quasi-Monte Carlo rules. 

In a number of practical applications one requires point sets with low star-discrepancy whose cardinality $N$ is not too large compared to $s$. This arises naturally, for instance, in estimating the expectation value of a linear functional of the solution to a partial differential equation with random coefficients \cite{KSS}. From \cite[Theorem~8]{KSS} one can see that in order to reduce the overall error, one needs to balance $N$ and $s$ and this balancing is of the form $N = s^\kappa$, for some $0 < \kappa \le 1$. Another case arises when the dimension $s$ is very large. For instance, in some applications from financial mathematics, the dimension $s$ can be several hundreds, see for instance \cite{PT}. If $s = 100$, then $2^s \approx 10^{30}$. Due to the limitations of the current technology, the number of points $N$ we can use is much smaller than $2^s$ in this case. In some instances of these applications, one can still achieve a convergence rate $N^{-\delta}$ with $\delta > 1/2$, but if the problem is 
more difficult the best we can get is $N^{-1/2}$. In these situations one would like to reduce the constant factor which depends on the dimension $s$, thus point sets whose star-discrepancy depends only weakly on the dimension $s$ are useful. At the end of the paper we discuss another situation where the dependence on the dimension is important.


\section{The inverse of the star-discrepancy problem}

We review the current literature on the inverse of the star-discrepancy problem as first studied in \cite{HNWW}. To analyze the problem systematically the so-called {\it inverse of the star-discrepancy} is defined as $$N(s,\varepsilon)=\min\{N \in \mathbb{N}\ : \ {\rm disc}(N,s) \le \varepsilon\}\ \ \mbox{ for  } \ s \in \mathbb{N} \mbox{ and } \varepsilon \in (0,1].$$ This is the minimal number of points which is required to achieve a star-discrepancy less than $\varepsilon$ in dimension $s$. The following theorem is the first classic result in this direction.

\begin{theorem}[Heinrich, Novak, Wasilkowski and Wo\'{z}niakowski \cite{HNWW}]\label{thm7.1}
We have  
\begin{equation}\label{mindihnww}
{\rm disc}(N,s) \ll \sqrt{\frac{s}{N}} \ \ \mbox{ for all }\ \ N,s \in \mathbb{N}.
\end{equation}
Hence $$N(s,\varepsilon) \ll s \varepsilon^{-2}\ \ \mbox{ for  all }\ s \in \mathbb{N}\ \mbox{ and }\ \varepsilon>0.$$ 
\end{theorem}

The bound \eqref{mindihnww} does not achieve the optimal rate of convergence for fixed dimension $s$ as the number of points $N$ goes to $\infty$. However, the dependence on the dimension $s$ is much weaker than in \eqref{star_bounds}. Thus such point sets are more suited for integration problems where the dimension $s$ is large.

The proof of Theorem~\ref{thm7.1} is based on the probabilistic method. It is shown that the probability, that the absolute local discrepancy $|\Delta(\boldsymbol{t})|$ of a randomly chosen point set is larger than a certain quantity $\delta$, is extremely small. Then one applies a union bound over all $\boldsymbol{t} \in [0,1]^s$ and chooses $\delta$ such that this union bound is strictly less then one, which then implies the result. In this particular instance the authors of \cite{HNWW} used a large deviation inequality for empirical processes on Vapnik-\v{C}ervonenkis classes due to Talagrand and Haussler. Details can be found in \cite{DHP,HNWW}. A simplified proof which leads in addition to explicit constants was given recently by Aistleitner \cite{aist2011}. 

It is also known that the dependence on the dimension $s$ of the inverse of the star-discrepancy cannot be improved. Hinrichs \cite{Hin2004} proved the existence of constants $c, \varepsilon_0 >0$ such that $$N(s,\varepsilon) \ge c s \varepsilon^{-1} \ \ \mbox{ for all } \varepsilon \in (0,\varepsilon_0) \mbox{ and } s \in \mathbb{N}$$ and ${\rm disc}(N,s) \ge \min(\varepsilon_0,c s/n)$. The exact dependence of $N(s,\varepsilon)$ on $\varepsilon^{-1}$ is still an open question which seems to be very difficult.

Doerr~\cite{Doerr} on the other hand showed that the star-discrepancy of a random point set is at least of order $\sqrt{s/N}$, which shows that the upper bound of \cite{HNWW} is asymptotically sharp. 

A similar but slightly weaker result compared to Theorem~\ref{thm7.1} is the following:

\begin{theorem}[Heinrich, Novak, Wasilkowski and Wo\'{z}niakowski \cite{HNWW}]\label{ditract}
We have 
\begin{equation}\label{nhwwrhs}
{\rm disc}(N,s) \ll  \sqrt{\frac{s}{N}} \ \sqrt{\log s+\log N}  \ \ \mbox{ for all }\ N,s \in \mathbb{N},
\end{equation}
and  
\begin{equation*}\label{hnwwmindisc1}
N(s,\varepsilon) \ll s \varepsilon^{-2} \log(s/\varepsilon) \ \ \mbox{ for  all }\ s \in \mathbb{N}\ \mbox{ and }\ \varepsilon>0.
\end{equation*}
\end{theorem} 

The proof of this result is based on similar ideas as used in the proof of the previous theorem, but instead of the result of Talagrand and Haussler, here the authors of \cite{HNWW} used Hoeffdings inequality, which is an estimate for the deviation from the mean for sums of independent random variables. We give a short sketch of the proof which offers some insights. More details can be found in \cite{DP10,HNWW,LP14}. 

\noindent{\it Sketch of the proof of Theorem~\ref{ditract}.} Hoeffdings inequality (in the form required here) states that if $X_0,\dots,X_{N-1}$ are independent real valued random variables with mean zero and $|X_i| \le 1$ for $i=0,\dots,N-1$ almost surely, then for all $t>0$ we have
 $$
   \mathrm{Prob} \left( \left|\sum_{i=0}^{N-1} X_i \right| > t \right) \le 2 \exp\left( - \frac{t^2}{2N}\right).
 $$

Now let $P_{N,s}=\{\boldsymbol{x}_0,\dots,\boldsymbol{x}_{N-1}\}$ where $\boldsymbol{x}_0,\dots,\boldsymbol{x}_{N-1}$ are independent and uniformly distributed in $[0,1)^s$. We want to show that
$$ \mathrm{Prob} \left( D_N^{\ast}(P_{N,s}) \le 2 \varepsilon \right) > 0 $$
where $2 \varepsilon$ is the right hand side in Theorem~\ref{nhwwrhs}.
That amounts to the task to show that the event
$$|\Delta(\boldsymbol{x})| > 2 \varepsilon \ \ \mbox{at least for one } \boldsymbol{x}\in[0,1)^s$$
has a probability smaller then 1. 
These are infinitely many constraints, but it can be shown that 
$|\Delta(\boldsymbol{x})| > 2 \varepsilon$ implies  $|\Delta(\boldsymbol{y})| > \varepsilon$ for one of the points   
in a rectangular equidistant grid of mesh size $1/m$ with $m=\lceil s/\varepsilon \rceil$.
 Actually, this holds either for the grid point directly below left or up right from $\boldsymbol{x}$.
 Since this grid has cardinality $(m+1)^s$, a union bound shows that it is enough to prove
 $$  \mathrm{Prob} \left(|\Delta(\boldsymbol{x})| >  \varepsilon \right) < (m+1)^{-s} $$
 for every $\boldsymbol{x} \in [0,1)^s$.
 But now  
 $$ N \Delta(\boldsymbol{x}) = \sum_{i=0}^{N-1} \left( {\mathbf 1}_{[\boldsymbol{0},\boldsymbol{x})}(\boldsymbol{x}_i) -  \, \lambda_s([\boldsymbol{0},\boldsymbol{x})) \right) $$
 is the sum of the $N$ random variables $X_i={\mathbf 1}_{[\boldsymbol{0},\boldsymbol{x})}(\boldsymbol{x}_i) -  \, \lambda_s([\boldsymbol{0},\boldsymbol{x}))$, which have mean 0 and
 obviously satisfy $|X_i|\le 1$.
 So we can apply Hoeffding's inequality and obtain
 $$ \mathrm{Prob} \left(|\Delta(\boldsymbol{x})| >  \varepsilon \right)  = \mathrm{Prob} \left( \left|\sum_{i=0}^{N-1} X_i \right| > N\varepsilon \right) \le 2 \exp \left( \frac{-N \varepsilon^2}{2}\right) < (m+1)^{-s}, $$
 where the last inequality is satisfied for the chosen values of the parameters. $\qed$ 

The results in Theorem~\ref{thm7.1} and \ref{ditract} are only existence results. Until now no {\it explicit} constructions of $N$-element point sets $P_{N,s}$ in $[0,1)^s$ for which $D_N^{\ast}(P_{N,s})$ satisfy \eqref{mindihnww} or \eqref{nhwwrhs} are known. A first constructive approach was given by Doerr, Gnewuch and Srivastav~\cite{dgs05}, which is further improved by Doerr and Gnewuch~\cite{doegne06}, Doerr, Gnewuch, and Wahlstr\"om~\cite{DGW10} and Gnewuch, Wahlstr\"om and Winzen~\cite{GWW12}. There, a deterministic algorithm is presented that constructs point sets $P_{N,s}$ in $[0,1)^s$ satisfying $$D_N^\ast(P_{N,s}) \ll \sqrt{\frac{s}{N}} \ \sqrt{\log(N+1)}$$ in run-time $O(s \log(sN)(\sigma N)^s)$, where $\sigma=\sigma(s)=O((\log s)^2/(s \log \log s)) \rightarrow 0$ as $s \rightarrow \infty$ and where the implied constants in the $O$-notations are independent of $s$ and $N$. However, this is by far too expensive to obtain point sets for high dimensional applications. A slight improvement for the 
run time is presented in Doerr, Gnewuch, Kritzer and Pillichshammer~\cite{dgkp08}, but this 
improvement has to be payed with by a worse dependence of the bound on the star-discrepancy on the dimension.

\section{The weighted star-discrepancy}

In the paper \cite{SW98}, Sloan and Wo\'zniakowski introduced the notion of weighted star-discrepancy and proved a ``weighted'' Koksma-Hlawka inequality. The idea is that in many applications some projections are more important than others and that this should also be reflected in the quality measure of the point set. 

We start with some basic notation: let $[s]=\{1,2,\ldots ,s\}$ denote the set of coordinate indices. Let $\boldsymbol{\gamma}=(\gamma_j)_{j \ge 1}$ be a sequence of nonnegative reals. For $\mathfrak{u} \subseteq [s]$ we write $\gamma_{\mathfrak{u}}=\prod_{j \in \mathfrak{u}}\gamma_{j}$, where the empty product is one by definition. The real number $\gamma_{\mathfrak{u}}$ is the ``weight'' corresponding to the group of variables given by $\mathfrak{u}$. Let $|\mathfrak{u}|$ be the cardinality of $\mathfrak{u}$. For a vector $\boldsymbol{z} \in [0,1]^s$ let $\boldsymbol{z}_{\mathfrak{u}}$ denote
the vector from $[0,1]^{|\mathfrak{u}|}$ containing the components of $\boldsymbol{z}$
whose indices are in $\mathfrak{u}$. By $(\boldsymbol{z}_{\mathfrak{u}},1)$ we mean the vector
$\boldsymbol{z}$ from $[0,1]^s$ with all components whose indices are not in
$\mathfrak{u}$ replaced by 1. 

For an $N$-element point set $P_{N,s}$ in $[0,1)^s$ and given weights
$\boldsymbol{\gamma}=(\gamma_j)_{j \ge 1}$, the {\it weighted star-discrepancy} $D_{N,\boldsymbol{\gamma}}^{\ast}$ is given by
\begin{eqnarray*}
D_{N,\boldsymbol{\gamma}}^{\ast}(P_{N,s}) =\sup_{\boldsymbol{z} \in [0,1]^s} \max_{\emptyset \not=\mathfrak{u} \subseteq [s]}
\gamma_{\mathfrak{u}} |\Delta(\boldsymbol{z}_{\mathfrak{u}},1)|.
\end{eqnarray*}
If $\gamma_j = 1$ for all $j \ge 1$, then the weighted star-discrepancy coincides with the classical star-discrepancy.  

Quite similar to the classical case, we define the {\it $N$th minimal weighted star-discrepancy} $${\rm disc}_{\boldsymbol{\gamma}}(N,s)=\inf_{P_{N,s}}  D_{N,\boldsymbol{\gamma}}^{\ast}(P_{N,s})$$ and the {\it inverse of the weighted star-discrepancy} $$N_{\boldsymbol{\gamma}}(s,\varepsilon)=\min\{N \in \mathbb{N}\ : \ {\rm disc}_{\boldsymbol{\gamma}}(N,s) \le \varepsilon\}.$$ 

Now we recall two notions of tractability. Tractability means that we control the dependence of the inverse of the weighted star-discrepancy on $s$ and $\varepsilon^{-1}$ and rule out the cases for which $N_{\boldsymbol{\gamma}}(s,\varepsilon)$ depends exponentially on s or on $\varepsilon^{-1}$.

\begin{itemize}
\item We say that the weighted star-discrepancy is {\it polynomially tractable}, if there exist nonnegative real numbers $\alpha$ and $\beta$ such that 
\begin{equation}\label{pt}
N_{\boldsymbol{\gamma}}(s,\varepsilon) \ll  s^{\beta} \varepsilon^{-\alpha} \ \ \mbox{ for all } s\in \mathbb{N} \mbox{ and } \varepsilon \in (0,1).
\end{equation}
The infima over all $\alpha,\beta>0$ such that \eqref{pt} holds are called the $\varepsilon$-exponent and the $s$-exponent, respectively, of polynomial tractability. 
\item We say that the weighted star-discrepancy is {\it strongly polynomially tractable}, if there exists a nonnegative real number $\alpha$ such that 
\begin{equation}\label{spt}
N_{\boldsymbol{\gamma}}(s,\varepsilon) \ll \varepsilon^{-\alpha}\ \ \mbox{ for all } s\in \mathbb{N} \mbox{ and } \varepsilon \in (0,1).
\end{equation}
The infimum over all $\alpha>0$ such that \eqref{spt} holds is called the $\varepsilon$-exponent of strong polynomial tractability. 
\end{itemize} 
In both cases the implied constant in the $\ll$ notation is independent of $s$ and $\varepsilon$.

We collect some known results for the weighted star-discrepancy. The first result is an extension of Theorem~\ref{ditract} to the weighted star-discrepancy.
\begin{theorem}[Hinrichs, Pillichshammer, Schmid \cite{HPS2008}]\label{thHPS}
We have 
\begin{equation}\label{bdx1}
{\rm disc}_{\boldsymbol{\gamma}}(N,s) \ll \frac{\sqrt{\log s}}{\sqrt{N}}  \max_{\emptyset \not=\mathfrak{u}
\subseteq [s]} \gamma_{\mathfrak{u}} \sqrt{|\mathfrak{u}|}.
\end{equation}
\end{theorem}
Note that the result holds for every choice of weights. It is a pure existence result. Under very mild conditions on the weights, Theorem~\ref{thHPS} implies polynomial tractability with $s$-exponent zero. See \cite{HPS2008} for details. A slightly improved and numerically explicit version of Theorem~\ref{thHPS} can be found in the recent paper of Aistleitner~\cite{Aist}.

Next we state the following result: 

\begin{theorem}[Dick, Leobacher, Pillichshammer \cite{DLP2006}]\label{thHPS2}
For every prime number $p$, every $m \in \mathbb{N}$ and for given weights $\boldsymbol{\gamma}=(\gamma_j)_{j \ge 1}$ with $\sum_j \gamma_j < \infty$ one can construct (component-by-component) a $p^m$-element point set $P_{p^m,s}$ in $[0,1)^s$ such that for every $\delta >0$ we have $$D_{p^m,\boldsymbol{\gamma}}^{\ast}(P_{p^m,s})\ll_{\boldsymbol{\gamma},\delta} \frac{1}{p^{m(1-\delta)}}.$$
\end{theorem}

Note that the point set $P_{p^m,s}$ from Theorem~\ref{thHPS2} depends on the choice of weights. The result implies that the weighted star-discrepancy is strongly polynomially tractable with $\varepsilon$-exponent equal to one, as long as the weights $\gamma_j$ are summable. See \cite{DLP2006,dnp06,DP10,HPS2008} for more details. 

The next result (which follows implicitly from \cite{wang2}) is about Niederreiter sequences in prime-power base $q$. For the definition of Niederreiter sequences we refer to \cite{DP10,niesiam}. 

\begin{theorem}[Wang \cite{wang2}]\label{thm_afterWang}
For the weighted star-discrepancy of the first $N$ elements $P_{N,s}$ of an $s$-dimensional Niederreiter sequence in prime-power base $q$ we have $$D_{N,\boldsymbol{\gamma}}^{\ast}(P_{N,s}) \le \frac{1}{N} \max_{\emptyset \not= \mathfrak{u} \subseteq [s]} \prod_{j \in \mathfrak{u}} \left[ \gamma_j (C \ j \log(j+q) \log(qN))\right],$$ with a suitable constant $C>0$.
\end{theorem}

One can easily deduce from Theorem~\ref{thm_afterWang} that the weighted star-discrepancy of the Niederreiter sequence can be bounded independently of the dimension whenever the weights satisfy $\sum_j \gamma_j j \log j < \infty$. This implies strong polynomial tractability with $\varepsilon$-exponent equal to one. A similar result can be shown for Sobol' sequences and for the Halton sequence (see \cite{wang1,wang2}). 

We also have the following recent existence result:
\begin{theorem}[Aistleitner \cite{Aist}]\label{thm_aist}
For product weights satisfying $\sum_{j} {\rm e}^{-c \gamma_j^{-2}} < \infty$, for some $c > 0$, we have
$${\rm disc}_{\boldsymbol{\gamma}}(N,s) \ll_{\boldsymbol{\gamma}} \frac{1}{\sqrt{N}}\ \ \mbox{ for all }\ s,N \in \mathbb{N}.$$ 
Consequently, the weighted star-discrepancy for such weights is strongly polynomially tractable, with $\varepsilon$-exponent
at most 2.
\end{theorem}

All results described so far have either been existence results of point sets with small star-discrepancy, or results for point sets with small star-discrepancy which can be obtained via computer search. The Ansatz via computer search remains difficult and is limited to a rather small number of points $N$ and dimensions $s$ (in fact, it is known that the computation of the star-discrepancy is $NP$-hard as shown by Gnewuch, Srivastav, and Winzen~\cite{GSW}, which makes it difficult to obtain good point sets via computer search). To make the random constructions useful in applications, Aistleitner and Hofer~\cite{AH14} show that with probability $\delta$ one can expect point sets with discrepancy of order $c(\delta) \sqrt{s/N}$. Another Ansatz for obtaining explicit constructions is contained in \cite{Tem}.

In the following section we discuss results for explicit constructions of point sets. 

\section{The weighted star-discrepancy of Korobov's $p$-sets}

Let $p$ be a prime number. For a nonnegative real number $x$ let $\{ x \} = x - \lfloor x \rfloor$ denote the fractional part of $x$. For vectors we use this operation component-wise. 

We consider the so-called {\it $p$-sets} in $[0,1)^s$, a term which goes back to Hua and Wang~\cite{huawang}: 
\begin{itemize}
\item Let $P_{p,s}=\{\boldsymbol{x}_0,\ldots,\boldsymbol{x}_{p-1}\}$ with $$\boldsymbol{x}_n=\left(\left\{\frac{n}{p}\right\},\left\{\frac{n^2}{p}\right\},\ldots,\left\{\frac{n^s}{p}\right\}\right)\ \ \ \mbox{ for }\ n=0,1,\ldots,p-1.$$ The point set $P_{p,s}$ was introduced by Korobov \cite{kor1963}.
\item Let $Q_{p^2,s}=\{\boldsymbol{x}_0,\ldots,\boldsymbol{x}_{p^2-1}\}$ with $$\boldsymbol{x}_n=\left(\left\{\frac{n}{p^2}\right\},\left\{\frac{n^2}{p^2}\right\},\ldots,\left\{\frac{n^s}{p^2}\right\}\right)\ \ \ \mbox{ for }\ n=0,1,\ldots,p^2-1.$$ The point set $Q_{p,s}$ was introduced by Korobov \cite{kor1957}.
\item Let $R_{p^2,s}=\{\boldsymbol{x}_{a,k}\ : \ a,k \in \{0,\ldots,p-1\}\}$ with $$\boldsymbol{x}_{a,k}=\left(\left\{\frac{k}{p}\right\},\left\{\frac{a k}{p}\right\},\ldots,\left\{\frac{a^{s-1} k}{p}\right\}\right)\ \ \ \mbox{ for }\ a,k=0,1,\ldots,p-1.$$ Note that $R_{p^2,s}$ is the multi-set union of all Korobov lattice point sets with modulus $p$. The point set $R_{p^2,s}$ was introduced by Hua and Wang (see \cite[Section 4.3]{huawang}).
\end{itemize}

\begin{theorem}[Dick and Pillichshammer \cite{DPpset}]\label{thm2pset}
Assume that the weights $\gamma_j$ are non-increasing.
\begin{enumerate}
\item If $\sum_j \gamma_j < \infty$,  then for all $\delta>0$ we have
\begin{align*} 
D_{p,\boldsymbol{\gamma}}^{\ast}(P_{p,s}) \ll_{\boldsymbol{\gamma},\delta} \frac{1}{p^{1/2 -\delta}}, \ \ D_{p^2,\boldsymbol{\gamma}}^{\ast}(Q_{p^2,s}) \ll_{\boldsymbol{\gamma},\delta} \frac{1}{p^{1 -\delta}}, \ \mbox{ and }
D_{p^2,\boldsymbol{\gamma}}^{\ast}(R_{p^2,s}) \ll_{\boldsymbol{\gamma},\delta} \frac{1}{p^{1 -\delta}},
\end{align*}
where in all cases the implied constant is independent of $p$ and $s$. This implies strong polynomial tractability.

\item If there exists a real $\tau > 0$ such that $\sum_j \gamma_j^{\tau} < \infty$, then for all $\delta>0$ we have
$$D_{p,\boldsymbol{\gamma}}^{\ast}(P_{p,s}) \ll_{\boldsymbol{\gamma},\delta} \frac{s}{p^{1/2 -\delta}},\ \
D_{p^2,\boldsymbol{\gamma}}^{\ast}(Q_{p^2,s}) \ll_{\boldsymbol{\gamma},\delta} \frac{s}{p^{1 -\delta}}, \ \mbox{ and }
D_{p^2,\boldsymbol{\gamma}}^{\ast}(R_{p^2,s}) \ll_{\boldsymbol{\gamma},\delta} \frac{s}{p^{1 -\delta}},$$
where in all cases the implied constant is independent of $p$ and $s$. This implies polynomial tractability.
\end{enumerate}
\end{theorem}

The proof of Theorem~\ref{thm2pset} is based on an Erd\H{o}s-Turan-Koksma-type inequality for the weighted star-discrepancy and the following estimates for exponential sums. For details we refer to \cite{DPpset}.
\begin{lemma}\label{le3}
Let $p$ be a prime number and let $s \in \mathbb{N}$. Then for all $h_1,\ldots,h_s\in \mathbb{Z}$ such that $p \nmid h_j$ for at least one $j \in [s]$ we have 
\begin{align}
\left|\sum_{n=0}^{p-1} \exp(2 \pi \mathtt{i}(h_1 n+h_2 n^2+\cdots+h_s n^s)/p)\right| \le & (s-1) \sqrt{p}, \label{psum1} \\ \left|\sum_{n=0}^{p^2-1} \exp(2 \pi \mathtt{i}(h_1 n+h_2 n^2+\cdots+h_s n^s)/p^2)\right| \le & (s-1) p, \ \mbox{ and } \nonumber \\ \left|\sum_{a=0}^{p-1} \sum_{k=0}^{p-1} \exp(2 \pi \mathtt{i} k (h_1 +h_2 a+\cdots+h_s a^{s-1})/p)\right| \le & (s-1) p. \nonumber
\end{align}
\end{lemma}

Inequality \eqref{psum1} is known as Weil bound \cite{weil} and is often used in the area of pseudo-random number generation. Constructions related to the $p$-sets have also been considered in \cite{D14}. All of these constructions are related to the generation of (streams of) pseudo-random numbers (rather than low-discrepancy point sets and sequences). This may not be so surprising since the original argument by Heinrich, et. al. \cite{HNWW} is based on random samples and pseudo-random numbers are designed to mimic randomness. We discuss pseudo-random number generators in the next section more generally.

\section{Complete uniform distribution and pseudo-random number generators}

Pseudo-random number generators are commonly used in computer simulations to replace real random numbers for various reasons. Those point sets are based on deterministic constructions with the aim to mimic randomness. A number of quality criteria are applied to such pseudo-random number generators to assess their quality. One such criterion is complete uniform distribution.

Let $u_1, u_2, \ldots \in [0,1]$ be a sequence of real numbers. For $s,N \in \mathbb{N}$ we define
\begin{equation*}
\boldsymbol{u}^{(s)}_n = (u_{(n-1)s+1}, \ldots, u_{ns}) \in [0,1]^s.
\end{equation*}
Then the sequence $(u_n)_{n\ge 1}$ is completely uniformly distributed if for every $s \ge 1$
\begin{equation*}
\lim_{N \to \infty} D^\ast ( \{\boldsymbol{u}^{(s)}_1, \ldots, \boldsymbol{u}^{(s)}_{N}\} ) = 0.
\end{equation*}
The concept of complete uniform distribution measures correlations between successive numbers $u_{i}, u_{i+s}, u_{i+2s}, \ldots$. Real random numbers are uncorrelated and thus their discrepancy goes to $0$ (in probability), and so one wants pseudo-random numbers with the same property. For instance, the classic construction by van der Corput $(\phi(n))_{n \ge 0}$ in base $2$, given by $$\phi(n) = \frac{n_0}{2} + \frac{n_1}{2^2} + \cdots + \frac{n_m}{2^{m+1}},$$ where $n$ has dyadic expansion $n = n_0 + n_1 2 + \cdots + n_m 2^m$, is not completely uniformly distributed, since $\phi(2n)$ lies in the interval $[0, 1/2)$, whereas $\phi(2n-1)$ lies in the interval $[1/2, 1)$.

Markov chain algorithms are a staple tool in statistics and the applied sciences for generating samples from distributions for which only partial information is available. As such they are an important class of algorithms which use pseudo-random number generators. In \cite{CDO} it was shown that if the random numbers which drive the Markov chain are completely uniformly distributed, then the Markov chain consistently samples the target distribution (i.e. yields the correct result). For instance, \cite[Theorem~4]{CDO} requires pseudo-random numbers $(u_n)_{n\ge 1}$ such that for every sequence of natural numbers $(s_N)_{N \ge 1}$ with $s_N = \mathcal{O}(\log N)$, we have
\begin{equation*}
\lim_{N\to\infty} D^\ast(\{\boldsymbol{u}^{(s_N)}_{1}, \ldots, \boldsymbol{u}^{(s_N)}_{N} \}) = 0.
\end{equation*}
In this case, bounds like \eqref{star_bounds} are not strong enough due to their dependence on the dimension. Even a bound of the form $C^s N^{-\delta}$, with $C > 1$ and some $\delta > 0$ which does not depend on the dimension $s$, is not strong enough, since for $s = c \log N$ with $c > \frac{\log C}{\delta}$, we have $C^s N^{-\delta} = C^{c \log N} N^{-\delta} = N^{-\delta + c \log C} \ge 1$ for all $N \in \mathbb{N}$ and so we do not get any convergence. 

Thus it would be interesting for applications to explicitly construct a deterministic sequence $(u_n)_{n \ge 1}$ such that, say
\begin{equation*}
D^\ast(\{ \boldsymbol{u}^{(s)}_1, \ldots, \boldsymbol{u}^{(s)}_N\}) \le C \frac{\sqrt{s \log N} }{\sqrt{N}} \quad \mbox{for  all } N, s \in \mathbb{N}.
\end{equation*}
The existence of such a sequence has already been shown in \cite[p. 684]{CDO} and an improvement has been shown in \cite{AW}. Such a sequence has good properties when viewed as a pseudo-random sequences but is also useful as a deterministic sequence in quasi-Monte Carlo integration.

\vspace{-0.3cm}

\section*{Acknowledgments}\label{sec:Acknowledgments}

J. Dick is support by a Queen Elizabeth II Fellowship from the Australian Research Council. F. Pillichshammer is supported by the Austrian Science Fund (FWF): Project F5509-N26, which is a part of the Special Research Program "Quasi-Monte Carlo Methods: Theory and Applications".

\end{document}